\documentclass[a4paper,10pt]{article}
\usepackage{t1enc}
\usepackage[latin2]{inputenc}
\usepackage{amsmath}
\usepackage{amssymb}
\usepackage{amsthm}
\usepackage{graphics}
\author{{\sc Bálint Vető}}
\title{\sffamily\bfseries The time evolution of permutations under random stirring}
\date{}
\theoremstyle{plain}
\newtheorem{thm}{\sffamily\bfseries Theorem}
\newtheorem{lemma}{\sffamily\bfseries Lemma}
\newtheorem{prop}{\sffamily\bfseries Proposition}

\renewcommand{\P}{\mathbb{P}}
\newcommand{\p}{\mathbb{P}_{\varepsilon,\delta,n}}
\renewcommand{\O}{\mathcal{O}}
\newcommand{\gnt}{\lfloor\sqrt nt\rfloor}
\renewcommand{\S}{{\bf S}}
\newcommand{\I}{1\!\!1}
\newcommand{\C}{{\bf C}}
\newcommand{\n}{^{(n)}}

\DeclareMathOperator{\Leb}{Leb}

\begin{document}
\maketitle

\begin{abstract}\noindent{\sffamily\bfseries Abstract.} We consider permutations
of $\{1,\dots,n\}$ obtained by $\gnt$ independent applications of random
stirring. In each step the same marked stirring element is transposed with
probability $1/n$ with any one of the $n$ elements. Normalizing by $\sqrt n$ we
describe the asymptotic distribution of the cycle structure of these
permutations, for all $t\ge0$, as $n\to\infty$.
\end{abstract}

\section{\sffamily\bfseries Introduction}
We consider the following random stirring mechanism: $n$ numbered balls are
given in the beginning on their corresponding numbered places. In each step,
independently, the first ball, which is referred to as the \emph{stirring
particle} or \emph{stirring element}, changes place with one of the $n$ balls
or stays unchanged with probability $1/n$. We investigate that permutation
which brings the balls from their initial place to their place after $i$ steps.

Formally, let $\pi\n(i)=T_i\n\circ T_{i-1}\n\circ\dots\circ T_1\n$ be a
permutation acting on the set $[n]:=\{1,\dots,n\}$. The permutations
$(T_i\n)_{i=1}^\infty$ are chosen independently with uniform distribution from
the $n-1$ transpositions moving the stirring particle and the identity
permutation.

Let $\sigma$ be a permutation of a finite set $S$, i.e.\ an $S\to S$ bijective
function. The cycles (orbits) of $\sigma$ are the sets of form
$\{v,\sigma(v),\sigma^2(v),\dots\}\subseteq S$ for some $v\in S$. The set $S$
is the disjoint union of its cycles. The cycle structure of $\sigma$ is the
sequence of the cardinalities of the different cycles in non-increasing order.

In our case one of the cycles can be distinguished from the others (namely the
cycle of the stirring element), which will be called the \emph{active cycle}.
For the total description it is enough to determine the distribution of the
cycle structure of the permutation $\pi\n(i)$ (regarding the active cycle
separately). This gives the distribution of the conjugacy class of $\pi\n(i)$
restricting ourself to the conjugation with permutations fixing the stirring
particle. The distribution of $\pi\n(i)$ is uniform within a fixed conjugacy
class.

\begin{figure}
\centering \resizebox{!}{65pt}{\includegraphics{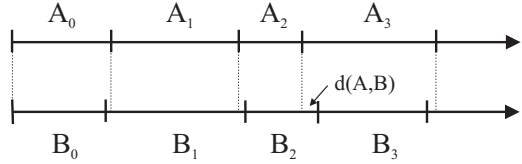}}

\caption{The metric on $\S$}\label{fig:d}
\end{figure}

We encode the permutation $\pi\n(i)$ with the vector
$\C\n(i):=(C_0\n(i),C_1\n(i),$ $C_2\n(i),\dots)$ where $C_0\n(i)$ denotes the
length of the active cycle, $C_1\n(i),C_2\n(i),\dots$ the lengths of those
cycles in non-increasing order which are already moved by one of the
transpositions $(T_j\n)_{j=1}^i$. Other $C_j\n(i)$-s are $0$.
$(\C\n(i))_{i=0}^\infty$ is a process on the state space
\begin{gather}
\S:=\{(s_0,s_1,s_2,\dots):s_n\in\mathbb{R},\quad s_n\ge0\quad n=0,1,2,\dots,\notag\\
s_1\ge s_2\ge\dots\ge s_n\ge\dots\mbox{ and }s_j>0\mbox{ for finitely many
$j$}\}
\end{gather}
with the distance
\begin{gather}
d(A,B):=\sup\left\{\left|\sum_{j=0}^k A_j-\sum_{j=0}^k
B_j\right|:k=0,1,2,\dots\right\}\label{def:d}
\end{gather}
where $A=(A_0,A_1,A_2,\dots)$ and $B=(B_0,B_1,B_2,\dots)$ are elements of \S.
(See Figure \ref{fig:d}.)
The ranking is not a natural part of the problem, but it facilitates studying
the model.

\begin{figure}
\centering \resizebox{!}{90pt}{\includegraphics{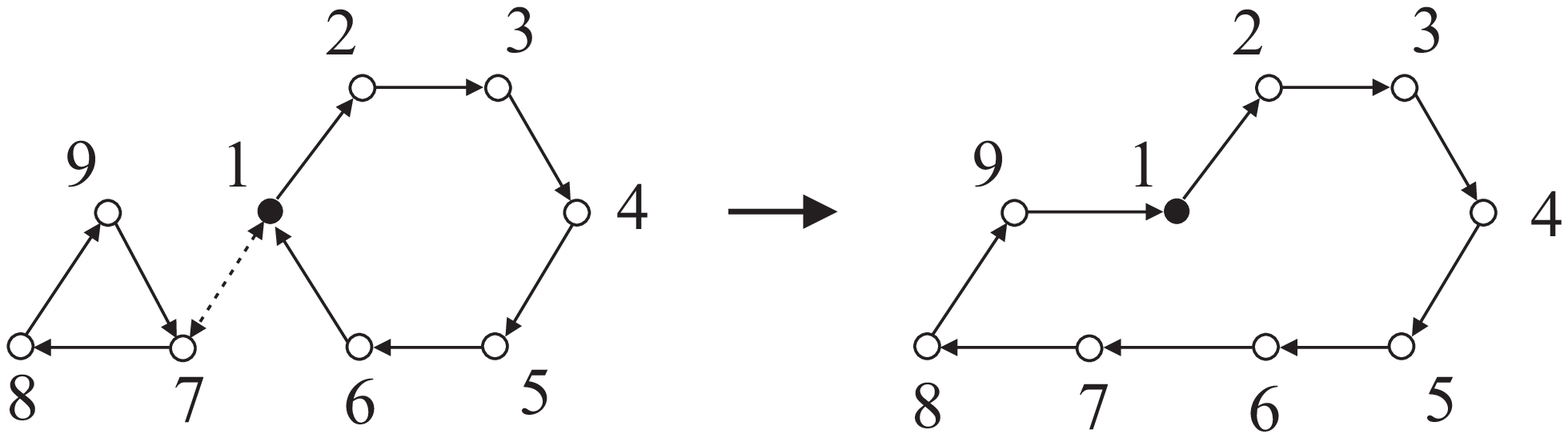}}
\resizebox{!}{120pt}{\includegraphics{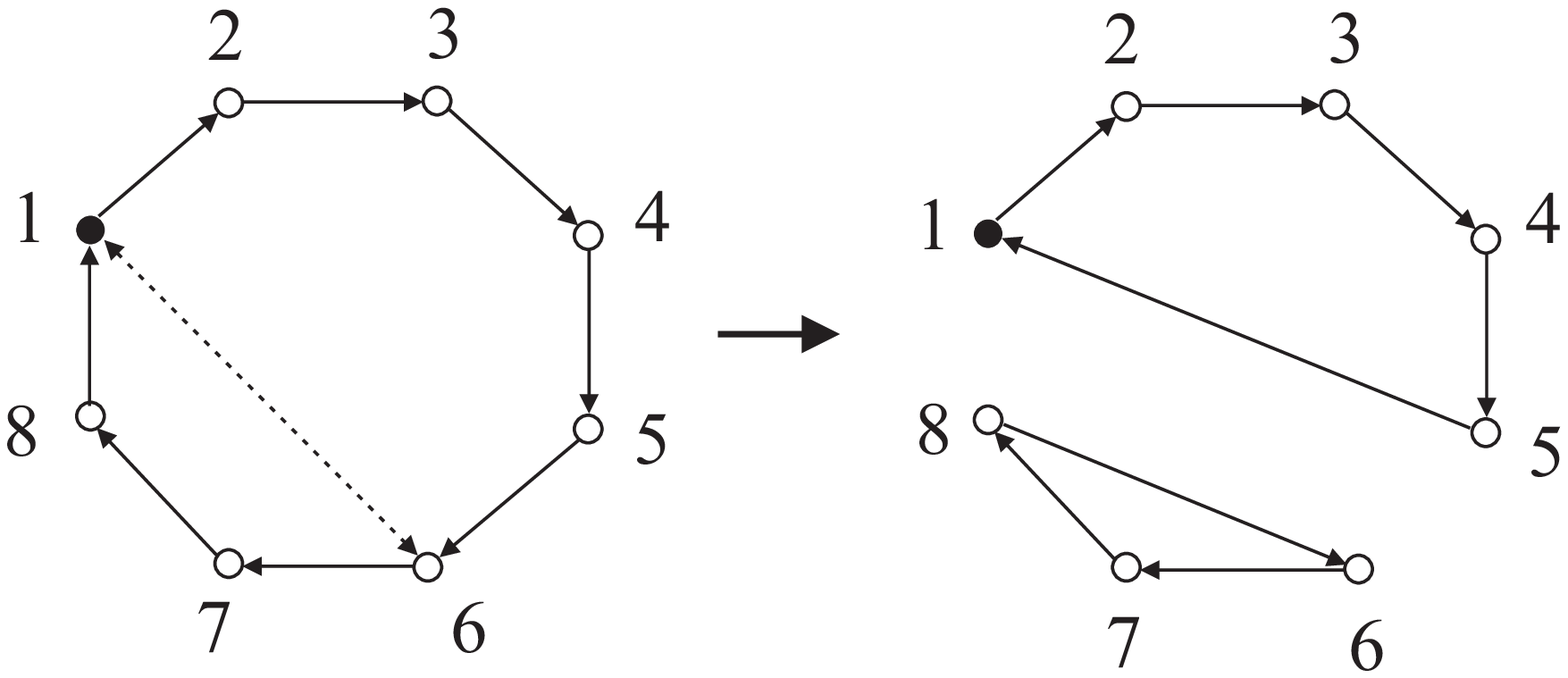}}\\[3ex]

\caption{Coagulation and fragmentation of cycles}\label{fig:cf1}
\end{figure}

At each step after applying a random transposition two types of changes may
happen in the cycle structure: merging of two distinct cycles or splitting of a
cycle in two. While different transpositions $(T_j\n)_{j=1}^i$ are applied
(meaning that the stirring particle chooses a new element in each step until
$i$), the cycle decomposition of $\pi\n(i)$ contains only fixed points and the
active cycle, which increases by one in each step: $\C\n(i)=(i+1,0,0,\dots)$.
If a transposition recurs, then the cycle splits in two, one of which will be
the new active cycle. If there are already more than one non-trivial cycles in
the decomposition, then the active cycle can merge another cycle. (See Figure
\ref{fig:cf1}.) The model realizes a coagulation-fragmentation process.

A reduction of the problem is to study the coagulation and fragmentation events
of the cycles together, because both of these events happen when the stirring
element steps to a place already visited. We investigate this simpler question
first. Then we introduce a continuous time process on \S, which turns out to be
the limit process. The convergence is proved by coupling. In Section 4.\ we
show that the stationary distribution of the underlying split-and-merge
transformation is the adequate modification of the Poisson\,--\,Dirichlet
distribution. (See the definition later.)

A similar model is studied by Schramm in \cite{schramm}. He chooses
$(T_i)_{i=1}^\infty$ to be independent random transpositions with uniform
distribution from all possible transpositions of the set $[n]$. The limit
distribution of the proportions of the giant cycles in the permutation
$\Pi(t)=T_t\circ\dots\circ T_1$ after $t=cn$ steps as $n\to\infty$ is
identified (where $c>1/2$ is a constant).

This result is in accordance with the classical theory of the random graphs
derived from Erdős \cite{erdos}. Let us consider the random graph $G(t)$ on the
vertex set $[n]$ where $\{u,v\}$ is an edge in it if and only if the
transposition $(u,v)$ appears in $\{T_1,\dots,T_t\}$. By the Erdős-Rényi
Theorem \cite{erdosrenyi} the graph $G(t)$ has a giant connected component only
in the case $t/n=c>1/2$ similarly to the condition on the random permutations.
(For random graphs and random graph processes see \cite{jlr} and
\cite{stepanov}).

In Schramm's paper the vector of the cycle sizes of $\Pi(t)$ in non-increasing
order normalized by the magnitude of the giant connected component of $G(t)$
converges in distribution to Poisson\,--\,Dirichlet distribution with parameter
$1$ after $t=cn$ steps ($c>1/2$) as $n\to\infty$. That is the limit
distribution of the relative cycle sizes in a random permutation chosen
uniformly from all permutations of $[n]$ as $n\to\infty$. Thus for large $n$
the permutation $\Pi(t)$ behaves on the giant connected component of the
Erdős-Rényi graph $G(t)$ as a uniform permutation.

Our paper is motivated by Tóth in connection with the quantum-physical
applications of the problem \cite{toth}. Angel analysed Tóth's random walk
model on regular trees in \cite{angel}. For similar random stirring models see
also \cite{durrett} and \cite{zeitouni}.

\section{\sffamily\bfseries Return times of the stirring particle} The movement of the
stirring particle is a random walk $(B_i\n)_{i=0}^\infty$ on the set $[n]$,
which is homogeneous in space and time. Let
\begin{gather}
V_i\n:=\#\{k:k\le i,\exists j<k:B_j\n=B_k\n\}\label{def:V}
\end{gather}
be the number of the \emph{return}s until the $i$th step to places already
visited by the random walk $(B_j\n)_{j=0}^\infty$. We also include those steps
when the stirring particle keeps its place.

After the $i$th step the stirring element has already visited exactly
$i+1-V_i\n$ places (including the starting point), so the transition
probabilities of the Markov-chain $(V_i\n)_{i=0}^\infty$ are
\begin{gather}
\P\left(V_{i+1}\n-V_i\n=1|V_i\n\right)=1-\P\left(V_{i+1}\n-V_i\n=0|V_i\n\right)
=\frac{i+1-V_i\n}n.\label{P:V-V}
\end{gather}

In order to get a non-trivial limit distribution the time of the processes
should be accelerated. As opposed to Schramm \cite{schramm}, in Theorem 1 the
scaling is $\sqrt n$. This means that we describe the beginning of the
evolution, because after $\O(\sqrt n)$ steps the bulk of the elements is still
unchanged. Simultaneously we
normalize the cycle sizes with $\sqrt n$ and we let $n\to\infty$.\\

From now on we investigate the limit of the vectors
$\left(\frac{\C\n(\gnt)}{\sqrt n}\right)_{t\ge0}$ as $n\to\infty$, where the
division is meant coordinatewise, namely $\frac{\C\n(\gnt)}{\sqrt
n}:=(\frac{C_0\n(\gnt)}{\sqrt n},$ $\frac{C_1\n(\gnt)}{\sqrt n},\dots)$.
Elementary calculations, similar to the classical birthday problem, give the
following limit distribution of the returns. For limit theorems related to
generalizations of the birthday problem see also \cite{camarri-pitman}.

\begin{prop}
Let $(V_t)_{t\ge0}$ be an inhomogeneous Poisson point process with intensity
$\rho(t)=t$. Then
\begin{gather}
(V_{\gnt}\n)_{t\ge0}\stackrel{\rm
d}{\Rightarrow}(V_t)_{t\ge0}\qquad(n\to\infty)\label{t:gyk}
\end{gather}
in terms of the finite dimensional marginal distributions.
\end{prop}

\section{\sffamily\bfseries Coupling}
Much more can be stated for the above model. Not only $(V_{\gnt}\n)_{t\ge0}$,
but the sequence of the processes $\left(\frac{\C\n(\gnt)}{\sqrt
n}\right)_{t\ge0}$ converges. Moreover by means of coupling a stronger type of
convergence is realized.

The limit process is a natural continuous extension of the discrete processes
$(\C\n(i))_{i=0}^\infty$. For large $n$ the active coordinate $C_0\n(i)$
increases in the bulk of the steps (when no split or merge occurs). In the
times of jumps of $(V_i\n)_{i=0}^\infty$ a split or a merge happens depending
on the proportions of the cycle sizes as follows. The probability of a split in
the $i$th step, conditionally given that the stirring particle returns to a
place already visited, is
\begin{gather}
\frac{C_0\n(i-1)}{\sum_{m=0}^\infty C_m\n(i-1)}.\label{P:sp}
\end{gather}
The conditional probability of the merge of the $j$th cycle and the active one
is
\begin{gather}
\frac{C_j\n(i-1)}{\sum_{m=0}^\infty C_m\n(i-1)}.\label{P:m}
\end{gather}

We define an $\S$ valued continuous time stochastic process
$\C(t)=(C_0(t),C_1(t),$ $C_2(t),\dots)$ with c\`adl\`ag paths, which imitates
the above process. It is built on a Poisson point process $(V_t)_{t\ge0}$ with
intensity $\rho(t)=t$. Similarly to the discrete processes
$(\C\n(i))_{i=0}^\infty$ at the times of jumps of $(V_t)_{t\ge0}$ a split or a
merge event occurs with probability proportional to the coordinates of $\C$.

The initial state is $\C(0):=(0,0,0,\dots)$. The evolution of the process is
the following: the coordinate $C_0(t)$ increases with constant speed $1$
between the jumps of $(V_t)_{t\ge0}$. Let $\tau_k$ be the $k$th time of jump of
$(V_t)_{t\ge0}$, in other words $V_{\tau_k}=k$ and
$V_{\tau_k-}=\lim_{\varepsilon\downarrow0}V_{\tau_k-\varepsilon}=k-1$. Let
$(U_k)_{k=1}^\infty$ be i.i.d.\ random variables with uniform distribution on
$[0,1]$ independent of $(V_t)_{t\ge0}$. One of the next two actions occurs at
time $\tau_k$.

\begin{enumerate}
\item \emph{Split:} If
\begin{gather}
U_k\le\frac{C_0(\tau_k-)}{\sum_{m=0}^\infty C_m(\tau_k-)},\label{def:sp}
\end{gather}
then let $C_0(\tau_k):=U_k\sum_{m=0}^\infty C_m(\tau_k-)$, and the sequence
$(C_m(\tau_k))_{m=1}^\infty$ will be the collection of
$(C_m(\tau_k-))_{m=1}^\infty$ and $C_0(\tau_k-)-U_k\sum_{m=0}^\infty
C_m(\tau_k-)$ rearranged in decreasing order.

\item \emph{Merge:} Otherwise a unique index $j\ge1$ can be chosen a.s.\ via
\begin{gather}
\frac{\sum_{m=0}^{j-1}C_m(\tau_k-)}{\sum_{m=0}^\infty C_m(\tau_k-)}<U_k
\le\frac{\sum_{m=0}^j C_m(\tau_k-)}{\sum_{m=0}^\infty C_m(\tau_k-)}.
\label{def:j}\end{gather}

\noindent Let $C_0(\tau_k):=C_0(\tau_k-)+C_j(\tau_k-)$, and
$C_m(\tau_k):=C_m(\tau_k-)$ if $1\le m<j$, and $C_m(\tau_k):=C_{m+1}(\tau_k-)$
if $m\ge j$ restoring the decreasing order.
\end{enumerate}

Observe that $\sum_{m=0}^\infty C_m(t)=t$, but we did not use it to simplify
the formulas \eqref{def:sp} and \eqref{def:j} in the above definition because
the analogous discrete assertion is not true, compare with \eqref{P:sp} and
\eqref{P:m}.

The main result of this paper is that the normalized discrete processes
converge in probability to $(\C(t))_{t\ge0}$ in the following uniform sense in
terms of the distance defined by \eqref{def:d}.

\begin{thm}
There exists a probability space $(\Omega,\mathcal{F},\P)$, on which the
discrete processes $(\C\n(i))_{i=0}^\infty\quad n=1,2,\dots$ and the continuous
time process $(\C(t))_{t\ge0}$ can be jointly realized so that if $T>0$ is
fixed and $f(n)$ is any function tending to infinity with $n$, then
\begin{gather}
\P\left(\sup_{0\le t\le T} d\left(\C(t),\frac{\C\n(\gnt)}{\sqrt n}\right)
<\frac{f(n)}{\sqrt n}\right)\to1 \quad\mbox{as }n\to\infty.
\end{gather}
\end{thm}

\subsection{\sffamily\bfseries The convergence of the return process}
Let $(\Omega,\mathcal{F},\P)$ be such a probability space where a Poisson point
process $(V_t)_{t\ge0}$ with intensity $\rho(t)=t$ and the i.i.d.\ random
variables $(U_k)_{k=1}^\infty$ and $(Z_i\n)_{i,n=1}^\infty$ with uniform
distribution on $[0,1]$ are given independently of each other.

We have constructed the process $(\C(t))_{t\ge0}$ from $(V_t)_{t\ge0}$ and
$(U_k)_{k=1}^\infty$ earlier. We first re-create the processes
$(V_i\n)_{i=0}^\infty$ with the appropriate distributions on the new
probability space $(\Omega,\mathcal{F},\P)$. The main idea of the construction
is that we observe the process $(V_t)_{t\ge0}$ in $\frac1{\sqrt n}$ long time
intervals.

Let $X_i\n:=\I\left(V_{\frac i{\sqrt n}}-V_{\frac{i-1}{\sqrt
n}}\ge1\right)\quad i=1,2\dots \quad n=1,2,\dots$ be the indicators of the
increase of the process $(V_t)_{t\ge0}$, which are Bernoulli random variables
with respective parameters
\begin{gather}
p_i\n=1-\exp\left(-\frac{2i-1}{2n}\right)=\frac
in+\O\left(\frac{i^2}{n^2}\right).\label{def:p}
\end{gather}

\noindent The required parameter for the increase of $V_i\n$ is
\begin{gather} q_i\n=\frac in- \frac{V_{i-1}\n}n.\label{def:q}\end{gather}

We define the values of $V_i\n$ for fixed $n$ with induction on $i$. Let
$V_0^{(n)}:=0\quad n=1,2,\dots$ and
\begin{gather}
Y_i\n:=X_i\n- \I\left(p_i\n>q_i\n\right) \I\left(X_i\n=1\right)
\I\left(Z_i\n>\frac{q_i\n}{p_i\n}\right)\notag\\
+\I\left(p_i\n<q_i\n\right)\I\left(X_i\n=0\right)
\I\left(Z_i\n<\frac{q_i\n-p_i\n}{1-p_i\n}\right).\label{def:Y}
\end{gather}
We define $V_i\n:=V_{i-1}\n+Y_i\n$.

It is easy to see that the distribution of the new $(V_i\n)_{i=0}^\infty$ is in
accordance with \eqref{P:V-V}. Later on we say that a \emph{correction} happens
if the products of the indicators in \eqref{def:Y} do not disappear. We will
see that the total probability that a correction ever occurs is small if $n$ is
large enough. This gives an alternative proof of Proposition 1.

\begin{lemma}
Let $T>0$ be fixed and denote $0=\tau_0,\tau_1,\dots,\tau_\kappa$ the random
times of jumps of the process $(V_t)_{0\le t\le T}$ and denote
$0=\tau_0\n,\tau_1\n,\dots,\tau_{\kappa\n}\n$ that of the discrete process
$(V_{\gnt}\n)_{0\le t\le T}$ defined above. Then for sufficiently large $n$
with probability close to $1$ the number of the jumps are equal:
$\kappa=\kappa\n$. Furthermore, there exists a bijection between the jumps of
the processes in such a way that
\begin{gather}
|\tau_k-\tau_k\n|\le\frac1{\sqrt n}\quad k=1,\dots,\kappa
\end{gather}
holds with large probability.
\end{lemma}

For technical convenience we introduce the following events for fixed
$\varepsilon,\delta>0$:
\begin{gather}
E_\varepsilon:=\{V_{\lfloor\sqrt n T\rfloor}\n\le K_\varepsilon\quad
n=N_\varepsilon,N_\varepsilon+1,\dots\},\label{def:E}
\end{gather}
where $K_\varepsilon$ is a sufficiently large constant and $N_\varepsilon$ is a
threshold satisfying $\P(E_\varepsilon)\ge1-\varepsilon$. It makes sense by
Proposition 1. Let
\begin{gather}
M_\delta:=\{\min_{k:\tau_k\le T}\{\tau_k-\tau_{k-1}\}>\delta\}\cap
\{V_T-V_{T-\delta}=0\}\label{def:M}
\end{gather}
where $\tau_k$ is the time of the $k$th jump of the process $(V_t)_{0\le t\le
T}$ and $\tau_0=0$. It is elementary that
$\lim_{\varepsilon\downarrow0}\P(E_\varepsilon)=\lim_{\delta\downarrow0}\P(M_\delta)=1$.\\

\noindent{\sffamily\bfseries Proof of Lemma 1: }By \eqref{def:M}, on the event
$M_\delta$ the increment of the process $(V_t)_{0\le t\le T}$ on any interval
$\left[\frac i{\sqrt n}, \frac{i+1}{\sqrt n}\right]$ does not exceed $1$ if
$n>\frac1{\delta^2}$, hence $V_{\frac i{\sqrt n}}-V_{\frac{i-1}{\sqrt
n}}=X_i\n$. Since $V_{\frac i{\sqrt n}}\n-V_{\frac{i-1}{\sqrt n}}\n=Y_i\n$, it
is enough to prove that
\begin{gather}
\P(\{\exists i\le\lfloor\sqrt n T\rfloor:X_i\n\neq Y_i\n\}\cap
E_\varepsilon\cap M_\delta)\to0\quad(n\to\infty)\label{XneqY}
\end{gather}
for all fixed $\varepsilon,\delta>0$.

On the event $E_\varepsilon$, $X_i\n=1$ can be true for at most $K_\varepsilon$
many indices $i$. So the probability of the correction in the cases
$p_i\n>q_i\n$ satisfies
\begin{gather}
1-\frac{q_i\n}{p_i\n}=\O\left(\frac1{\sqrt n}\right)\quad(n\to\infty)
\end{gather}
using the power series of the exponential function and the equations
\eqref{def:p} and \eqref{def:q} estimating $p_i\n$ and $q_i\n$. If we add this
at most $K_\varepsilon$ many times, then the sum still goes to $0$ as
$n\to\infty$. A similar calculation shows that for an $i$, for which
$p_i\n<q_i\n$ holds, the probability of the correction is at most
\begin{gather}
\frac{q_i\n-p_i\n}{1-p_i\n}=\O\left(\frac1n\right)\quad(n\to\infty).
\end{gather}
Summing up for $i=1,\dots,\lfloor\sqrt n T\rfloor$ the total probability still
tends to $0$, as required.

\subsection{\sffamily\bfseries Splits and merges}
With the processes $(V_i\n)_{i=0}^\infty$ we have determined when a split or a
merge occurs, our task is now to define how it should happen. Similarly to the
definition of the limit process $(\C(t))_{t\ge0}$ we can prescribe the
evolution of the discrete processes $(\C\n(i))_{i=0}^\infty$ with the use of
the same independent uniform random variables $(U_k)_{k=1}^\infty$ as follows.
Let $C_0\n(0):=1,\quad C_m\n(0):=0\quad m=1,2,\dots$. The evolution of the
process $\C\n$ in the steps $i=1,2,\dots$ is described below:
\begin{itemize}
\item if $V_i\n-V_{i-1}\n=0$, then $C_0\n(i):=C_0\n(i-1)+1$ and other
coordinates unchanged,
\item if $V_i\n-V_{i-1}\n=1$ and $V_i\n=k$, then the uniform random variable
$U_k$ determines a unique index $j$ with probability $1$ as in \eqref{def:j}
via
\begin{gather}
\frac{\sum_{m=0}^{j-1}C_m\n(i-1)}{\sum_{m=0}^\infty C_m\n(i-1)}<U_k\le
\frac{\sum_{m=0}^j C_m\n(i-1)}{\sum_{m=0}^\infty C_m\n(i-1)}.\label{def:j\n}
\end{gather}
\end{itemize}
Similarly to the definition of the limit process
\begin{enumerate}
\item $j=0$: \emph{split}. If $U_k\sum_{m=0}^\infty C_m\n(i-1)<1$, then let
everything be unchanged: $\C\n(i):=\C\n(i-1)$, let us call this case
\emph{fictive split} (corresponding to the event that the stirring particle
keeps its place). Otherwise $C_0\n(i):=\lfloor U_k\sum_{m=0}^\infty
C_m\n(i-1)\rfloor$, let the broken fragment $C_0\n(i-1)-\lfloor
U_k\sum_{m=0}^\infty C_m\n(i-1)\rfloor$ add to the collection of nonactive
pieces $(C_m\n(i-1))_{m=1}^\infty$ to form the new ranked sequence
$(C_m\n(i))_{m=1}^\infty$.
\item $j>0$: \emph{merge}. Let $C_0\n(i):=C_0\n(i-1)+C_j\n(i-1)$ and for the
re-ranking $C_m\n(i):=C_m\n(i-1)$ if $0<m<j$, and $C_m\n(i):=C_{m+1}\n(i-1)$ if
$m\ge j$.
\end{enumerate}

It is easy to show that this new definition of $(\C\n(i))_{i=0}^\infty$
provides the same distribution as in the model generated by transpositions, so
we prove the convergence for these processes.\\

\noindent{\sffamily\bfseries Proof of Theorem 1: }Let $\varepsilon,\delta>0$ be
fixed. Let $A_n$ denote the event that the assertion of Lemma 1 holds for
$(V_{\gnt}\n)_{0\le t\le T}$. We restrict ourselves to the events
$E_\varepsilon\cap M_\delta\cap A_n$. Let us define a measure (which is not a
probability measure) on the sets $B\in\mathcal{F}$:
\begin{gather} \p(B):=\P(B\cap E_\varepsilon\cap M_\delta\cap A_n). \end{gather}
By Lemma 1 it is enough to show that for fixed $\varepsilon,\delta>0$ the
processes $\C(t)$ and $\frac{\C\n(\gnt)}{\sqrt n}$ are sufficiently close to
each other for large $n$ except a set with $\p$-measure tending to $0$ as
$n\to\infty$. The proof consists of the following steps:
\begin{enumerate}
\item We estimate the increase of the distance between $\C(t)$ and
$\frac{\C\n(\gnt)}{\sqrt n}$ between two successive split or merge events.
\item We introduce those events when the distance under discussion cannot be
estimated: the awkward events (defined later) and the fictive splits. We show
that they have small probability.
\item On the complementer event, which has probability tending to $1$ as
$n\to\infty$, we show that a merge does not increase the distance between
$\C(t)$ and $\frac{\C\n(\gnt)}{\sqrt n}$ very much.
\item We do this also for the splits.
\item We summarize the estimates.
\end{enumerate}

{\sc Step 1.} Let
\begin{gather}
d_k^-:=d\left(\C(\tau_k-),\frac{\C\n(\lfloor\sqrt n\tau_k\n\rfloor-1)}{\sqrt
n}\right),\quad d_k^+:=d\left(\C(\tau_k),\frac{\C\n(\lfloor\sqrt
n\tau_k\n\rfloor)}{\sqrt n}\right)\label{def:di}
\end{gather}
denote the distance between the discrete and continuous processes before and
after the time of the $k$th split or merge. (Recall that $\tau_k$ is the time
of the $k$th jump of $(V_t)_{0\le t\le T}$ and $\tau_k\n$ is that of
$(V_{\gnt}\n)_{0\le t\le T}$, which are close $\p$-almost surely by Lemma 1.)

While no split or merge occurs, the distance between the processes does not
increase very much. From Lemma 1 the difference between $\tau_k$ and $\tau_k\n$
can be at most $\frac1{\sqrt n}$. The discrete processes
$(\frac{\C\n(\gnt)}{\sqrt n})_{t\ge0}$ change only in the times which are
multiples of $\frac1{\sqrt n}$. Thus $\p$-almost surely
\begin{gather}d_k^-\le d_{k-1}^++\frac2{\sqrt n}.\label{ing}\end{gather}

{\sc Step 2.} From now on we investigate only the split or merge points of the
processes. At the $k$th time of jump of $(V_t)_{0\le t\le T}$ and
$(V_{\gnt}\n)_{0\le t\le T}$ we choose with the help of $U_k$ one of the
components of $\C(\tau_k-)$ and $\C\n(\lfloor\sqrt n\tau_k\n\rfloor-1)$ via
\eqref{def:j} and \eqref{def:j\n}. Let us call the possibility that these
components are of different indices an \emph{awkward event}. If an awkward
event or a fictive split (meaning that $U_k\sum_{m=0}^\infty C_m\n(i-1)<1$)
occurs, then we cannot estimate $d(\C,\C\n/\sqrt n)$. We will see that these
events have probability tending to $0$ as $n\to\infty$.

We can choose the components of $\C$ and $\C\n$ as follows. We set the
coordinates of the vector $\C(\tau_k-)/\sum_{m=0}^\infty C_m(\tau_k-)$ to the
real line from the origin one after another, which gives a partition of the
unit interval $[0,1]$. We do this also with the coordinates of
\begin{gather}
\frac{\C\n(\lfloor\sqrt n\tau_k\n\rfloor-1)} {\sum_{m=0}^\infty
C_m\n(\lfloor\sqrt n\tau_k\n\rfloor-1)}.\label{C/sumC}
\end{gather}

\noindent Let $W_k$ denote the set of those points in $[0,1]$ which are covered
by the coordinates of $\C$ and $\C\n$ of different indices. The probability of
the awkward events (which is an upper estimate for their $\p$-measure) is
exactly the Lebesgue measure of $W_k$.

We know that $\sum_{m=0}^\infty C_m(t)=t$ for all $t\ge0$. From the
construction
\begin{gather}
\frac{\gnt-K_\varepsilon}{\sqrt n}\le\frac{\sum_{m=0}^\infty C_m\n(\gnt)}{\sqrt
n} \le t\quad\mbox{if}\quad t\in[0,T],\label{becs:sumC}
\end{gather}
because at the split or merge points (occurring at most $K_\varepsilon$ many
times) the total length of the discrete process does not increase. From this
\begin{gather}
\left|\sum_m C_m(\tau_k-)-\frac{\sum_m C_m\n(\lfloor\sqrt n\tau_k\n\rfloor-1)}
{\sqrt n}\right|\hspace{10em}\notag\\\hspace{5em}
\le\left|\tau_k-\frac{\lfloor\sqrt{n}\tau_k\n\rfloor-1-K_\varepsilon}{\sqrt
n}\right| \le\frac{K_\varepsilon}{\sqrt n}+\frac3{\sqrt n} \label{becs:C-C}
\end{gather}
follows using Lemma 1.

From the above it is an elementary exercise to show that the distance between
the corresponding dividing points of the partitions of $[0,1]$ generated by
$\C(\tau_k-)/\tau_k$ and by the vector \eqref{C/sumC} can be, respectively, at
most
\[\frac{d_k^-+\frac{K_\varepsilon+3}{\sqrt n}}{\tau_k}.\]

Since the number of coordinates is at most $K_\varepsilon$, this provides the
following upper bound:
\begin{gather}
\Leb(W_k)\le \frac{d_k^-+\frac{K_\varepsilon+3}{\sqrt n}}{\tau_k}
K_\varepsilon\le \frac{d_k^-+\frac{K_\varepsilon+3}{\sqrt n}}{\delta}
K_\varepsilon,
\end{gather}
where we used the fact that $\tau_k=\sum_{m=0}^\infty C_m(\tau_k-)\ge\delta$
holds for $k=1,2,\dots$ on the event $M_\delta$. This yields
\begin{gather}
\p(\mbox{awkward event at $\tau_k$})\le\Leb(W_k)\le
\frac{d_k^-+\frac{K_\varepsilon+3}{\sqrt n}}\delta
K_\varepsilon.\label{P:rossz}
\end{gather}

Furthermore
\begin{gather}
\p(\mbox{fictive split at $\tau_k\n$})\le\frac{\frac1{\sqrt
n}}{\sum_{m=0}^\infty C_m\n(\lfloor\sqrt n\tau_k\n\rfloor-1)/\sqrt
n}\le\frac2{\delta\sqrt n},\label{P:fikt}
\end{gather}
if $n$ is large enough by \eqref{becs:sumC}. So we conclude that
\begin{gather}
\p(\mbox{awkward event or fictive split at the $k$th split or merge
point})\hspace{5em}\notag\\\hspace{10em}\le \frac{K_\varepsilon}\delta d_k^-
+\frac{K_\varepsilon^2+3K_\varepsilon+2}{\delta\sqrt n}.\label{P:ke}
\end{gather}

\begin{figure}
\centering \resizebox{!}{60pt}{\includegraphics{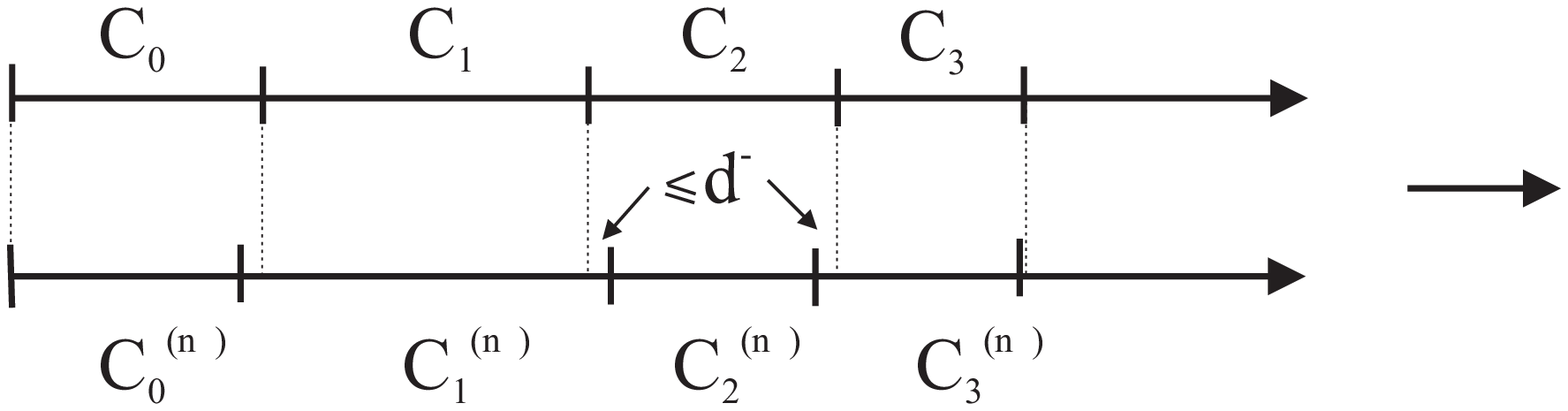}\hspace{2em}
\includegraphics{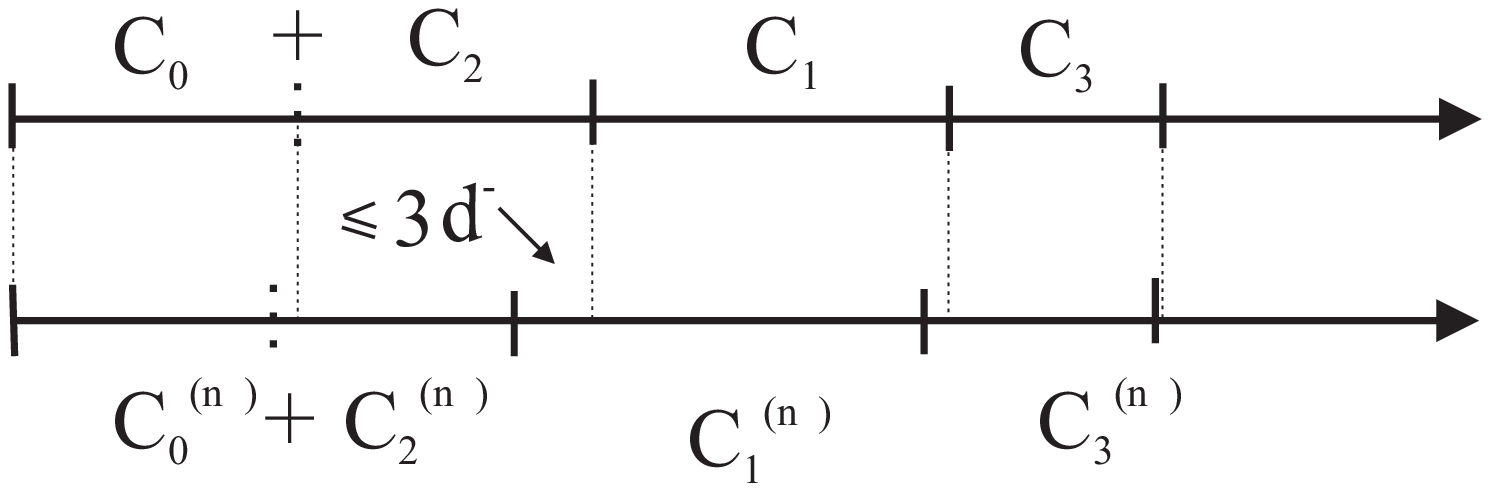}}
\caption{The piece $C_2$ merges $C_0$ parallel with the $C_2\n-C_0\n$
coagulation}\label{fig:coag}
\end{figure}

{\sc Step 3.} In the case when the random variable $U_k$ chooses the same
components of $\C$ and $\C\n$ and it is not the active coordinate, i.e.\ there
is a merge in both processes (see Figure \ref{fig:coag}), then
\begin{gather} d_k^+\le3d_k^-.\label{coag}\end{gather}

\begin{figure}
\centering \resizebox{!}{80pt}{\includegraphics{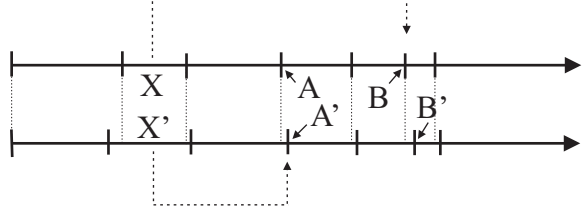}}

\caption{Split: the broken pieces from the coordinate $0$ are $X$ and $X'$
which have to be moved to places $A'$ and $B$}\label{fig:frag}
\end{figure}

{\sc Step 4.} If a (non-fictive) split occurs in the discrete and continuous
processes, then using inequality \eqref{becs:C-C} we have
\begin{gather}
\left|C_0(\tau_k)-\frac{C_0\n(\lfloor\sqrt{n}\tau_k\n\rfloor)}{\sqrt
n}\right|\notag\\
\le\left|U_k\sum_{m=0}^\infty C_m(\tau_k-)-\frac{\lfloor U_k\sum_{m=0}^\infty
C_m\n(\lfloor\sqrt n\tau_k\n\rfloor-1)\rfloor}{\sqrt n}\right|
\le\frac{K_\varepsilon+3}{\sqrt n}+\frac1{\sqrt n}.
\end{gather}
This is why the broken pieces from $C_0(\tau_k-)$ and $C_0\n(\lfloor\sqrt
n\tau_k\n\rfloor-1)/\sqrt n$ (denoted by $X$ and $X'$ on Figure \ref{fig:frag})
can differ at most $d_k^-+\frac{K_\varepsilon+4}{\sqrt n}$: the difference can
be $\frac{K_\varepsilon+4}{\sqrt n}$ between the left end points and at most
$d_k^-$ between the right end points.

It is possible that the two broken pieces do not come to the same place in the
decreasing order of the coordinates. This case is shown on Figure
\ref{fig:frag}. Then we move first both $X$ and $X'$ to the closer of the final
places of them in the decreasing order (to the places $A$ and $A'$ on the
figure). Because $|A-A'|\le d_k^-$, the result is two vectors (the
modifications of $\C(\tau_k-)$ and $\C\n(\lfloor\sqrt n\tau_k\n\rfloor-1)/\sqrt
n$, but one of them is not necessarily in decreasing order), which have
$d(\cdot,\cdot)$-distance at most $2d_k^-+\frac{K_\varepsilon+4}{\sqrt n}$ more
then before this modification.

In the second step we move $X$ from $A$ to $B$ (see Figure \ref{fig:frag}). The
lengths of the parts between $A$ and $B$ are at least $|X|$ and at most
$|X'|+2d_k^-\le |X|+3d_k^-+\frac{K_\varepsilon+4}{\sqrt n}$. So any two of
these parts have lengths differing at most $3d_k^-+\frac{K_\varepsilon+4}{\sqrt
n}$. Swapping $X$ always with its right neighbour until hitting place $B$, the
number of the swaps is at most $K_\varepsilon$, and at each swap the distance
can increase at most $3d_k^-+\frac{K_\varepsilon+4}{\sqrt n}$, so we have
\begin{gather}
d_k^+\le 2d_k^-+\frac{K_\varepsilon+4}{\sqrt n}+
K_\varepsilon\left(3d_k^-+\frac{K_\varepsilon+4}{\sqrt n}\right).\label{frag}
\end{gather}

{\sc Step 5.} Summing up the estimates \eqref{ing}, \eqref{coag} and
\eqref{frag} we get easily the following recursive bound:
\begin{gather}
d_k^+\le\max(3,2+3K_\varepsilon)d_{k-1}^+
+\frac{K_\varepsilon^2+5K_\varepsilon+4+2\max(3,2+3K_\varepsilon)}{\sqrt
n}=:ad_{k-1}^++\frac b{\sqrt n}.
\end{gather}
Hence $\sup_{0\le k\le K_\varepsilon}d_k^{\pm}
\le\left(\sum_{i=0}^{K_\varepsilon}ba^i\right)\frac1{\sqrt n}$. Considering the
results of steps 1 and 2 the assertion of the theorem follows.

\section{\sffamily\bfseries Stationary distribution and generalizations}\label{open}
It is a natural question to identify the stationary distribution of our
stirring process. This means that we look at the asymptotic behaviour of the
process $(\C\n(\lfloor nt\rfloor)/n)_{t\ge0}$. Observe that the time scale is
of order $n$, i.e.\ the time scale when the stirring element has already
visited the bulk of the $n$ places. This setup is the same as that of the
problem studied by Schramm in \cite{schramm}, but different from the phenomenon
described by Theorem 1.

In this section we consider the following split-and-merge transformation
corresponding to the stirring generated by random transpositions. Let
$\C=(C_0,C_1,$ $C_2,\dots)\in\S$ be a random probability distribution, i.e.\
$\sum_m C_m=1$ almost surely. $C_0$ is the active component. Let $U$ be a
random variable with uniform distribution of $[0,1]$ which is independent of
$\C$. If $U\le C_0$, then the $C_0$ splits, i.e.\ the new active component will
be $U$ and $(C_0-U,C_1,C_2,\dots)$ will be the remaining components after
restoring the decreasing order. If $\sum_{m=0}^{j-1}C_m<U\le\sum_{m=0}^j C_m$,
then $C_0$ merges with $C_j$ similarly to (\ref{P:sp}-\ref{def:j}) because
$\sum_m C_m=1$.

In limit theorems of random partitions and permutations the following
distribution appears often. Let the random variables $W_1,W_2,\dots$ be
independent with uniform distribution on $[0,1]$. Let $(Q_1,Q_2,\dots)$ be the
decreasing rearrangement of the random variables
\[(P_1,P_2,\dots):=(W_1,(1-W_1)W_2,(1-W_1)(1-W_2)W_3,\dots).\]
Then the random sequence $(P_1,P_2,\dots)$ has GEM($1$) distribution after
Griffiths, Engen and McCloskey. $(Q_1,Q_2,\dots)$ has
\emph{Poisson\,--\,Dirichlet} distribution with parameter $1$, abbreviated
PD($1$). For more about this family of distributions see \cite{pitman}.

Let $(p_1,p_2,\dots)$ be a random probability distribution. We construct its
size biased permutation. Let $U_1,U_2,\dots$ be i.i.d.\ uniform random
variables on $[0,1]$ independently of $(p_1,p_2,\dots)$. Let $I_j$ be the
unique index for which $\sum_{i=1}^{I_j-1}p_i\le U_j<\sum_{i=1}^{I_j}p_i$. Let
$J_k$ denote the $k$th smallest integer $m$ satisfying
$I_m\notin\{I_1,I_2,\dots,I_{m-1}\}$. Then the vector $(p_{J_1},p_{J_2},\dots)$
is called the size biased permutation of $(p_1,p_2,\dots)$. It is well known
that the size biased permutation of a random partition with PD($1$)
distribution has GEM($1$) distribution. See also \cite{pitman2}.

Consider the following probability distribution on \S. Let $(Q_1,Q_2,\dots)$
have PD($1$) distribution. Let $C_0$ be a size biased part from
$(Q_1,Q_2,\dots)$ (i.e.\ the first component of the size biased permutation of
$(Q_1,Q_2,\dots)$) corresponding to the active cycle and the rest
$(C_1,C_2,\dots)$ is the vector of the remaining $Q_j$-s in non-increasing
order. We denote by $\mu$ the distribution of $\C=(C_0,C_1,C_2,\dots)$.

\begin{thm}
The distribution $\mu$ is invariant under the above split-and-merge
transformation.
\end{thm}

\noindent{\sffamily\bfseries Proof: }By definition a random partition $\C$ with
distribution $\mu$ can be considered as follows. Let $W_1,W_2,\dots$ be i.i.d.\
uniform random variables on $[0,1]$ as in the definition of PD($1$). Because
the size biased permutation of PD($1$) is GEM($1$), we can suppose that for the
active component $C_0=W_1$ holds and $(C_1,C_2,\dots)$ is the decreasing
rearrangement of $((1-W_1)W_2,(1-W_1)(1-W_2)W_3,\dots)$. Let $\nu$ be the
distribution of the random partition obtained by the application of a stirring
step to $\C$.

If $U<W_1$ for the $[0,1]$-uniform random variable $U$, then the new non-active
components are $(W_1-U,(1-W_1)W_2,(1-W_1)(1-W_2)W_3,\dots)$ in decreasing
order. Conditionally on $\{U<W_1\}$ and on $U$, the variable $W_1$ is uniform
on $[U,1]$, thus the vector of the non-active components has PD($1$)
distribution scaled by $(1-U)$. It yields that $\nu$ conditioned on $\{U<W_1\}$
and on $U$ is the same as $\mu$ conditioned on the active component having size
$U$.

If $U>W_1$, then a coagulation occurs. Conditioned on $\{U>W_1\}$ and on the
value of $W_1$, the size of the component which merges $C_0$ has uniform
distribution on $[0,1-W_1]$, because it is a size biased component. We get the
same distribution, if we choose this component merging $C_0$ to be of length
$U-W_1$. Conditionally on $\{U>W_1\}$ and on $U$ the rest has PD($1$)
distribution scaled by $(1-U)$. Thus, a sample from $\nu$ conditioned on
$\{U>W_1\}$ and on $U$ has an active coordinate of size $U$ and the remaining
components with a scaled PD($1$) distribution.

Hence, a vector with distribution $\nu$ can be obtained by sampling $U$
uniformly on $[0,1]$, taking the active coordinate of length $U$ and taking a
scaled PD($1$) distribution on the rest. It shows that $\nu=\mu$, as required.\\

Theorem 2 proves that $\mu$ is a stationary measure for our process, but it is
not at all clear if this is the \emph{unique} stationary measure. The proof of
this would be the analogue of Schramm's result in \cite{schramm}.

A possible generalization of the model studied in this paper is the
\emph{multiple stirring}. It means that we consider more than one stirring
particles. For a fixed number $k$ of stirring elements an analogous limit
theorem can be proved with a coupling similarly to Theorem 1. The case, if the
number of the stirring elements depends on the size of the set $[n]$, might
also be worth studying (for example with $k(n)=n^\alpha$ where $0<\alpha<1$).
Of course, we need different scaling of time and space in this case.

An open question is for our original model to establish after how much time a
permutation can be regarded as a random permutation chosen with uniform
distribution, if it can be regarded at all. The solution of the problem in this
simply describable model is not obvious in the least. For more about this
problem in similar models see \cite{diaconis}.\\

\noindent{\sffamily\bfseries Acknowledgement. }I thank Bálint Tóth and Benedek
Valkó for initiating these investigations and for their permanent support
and useful comments while writing this paper. I am grateful to the referee for
pointing out the proof of Theorem 2, and I thank to Sándor Csörgő for his
helpful remarks on this paper.

\ \\\noindent{\sc B. Vető}, Institute of Mathematics, Technical University
Budapest, Egry J.\ u.\ 1, 1111 Budapest, Hungary; {\it e-mail:
}vetob@math.bme.hu
\end{document}